\documentclass[12pt]{article}
\usepackage{amssymb}
\usepackage{amsfonts}
\usepackage{amsmath}
\usepackage[usenames]{color}
\usepackage{mathrsfs}
\usepackage{amsfonts}
\usepackage{amssymb,amsmath}
\usepackage{CJK}
\usepackage{cite}
\usepackage{cases}
\usepackage{amsthm}\usepackage{multirow}
\usepackage{graphicx}
\usepackage{float}
\usepackage{makecell}
\usepackage{latexsym}
\usepackage{CJK}
\usepackage{bm}

\def\cl{\centerline}

\def\vs{\vspace*}

\def\ni{\noindent}

\numberwithin{equation}{section}
\newtheorem{theo}{Theorem}[section]

\newtheorem{coro}[theo]{Corollary}
\newtheorem{lemm}[theo]{Lemma}
\newtheorem{exam}[theo]{Example}

\newtheorem{rema}[theo]{Remark}

\newtheorem{remark}[theo]{Remark}
\newtheorem*{thmA}{Theorem A}

\begin{document}
\begin{center}
\cl{\large\bf \vs{6pt} Transnormal Functions and Focal Varieties on Finsler Manifolds \,{$^*\,$}}
\footnote {$^*\,$ Project supported by AHNSF (No.2108085MA11).
\\\indent\ \ $^\dag\,$ chenylwuhu@qq.com
}
\cl{Yali Chen$^1$$^\dag\,$, Qun He$^1$}

\cl{\small 1 School of Mathematical Sciences, Tongji University, Shanghai,
200092, China.}
\end{center}

{\small
\parskip .005 truein
\baselineskip 3pt \lineskip 3pt

\noindent{{\bf Abstract:}
In this paper, we study transnormal functions and their level sets and focal varieties on complete Finsler manifolds. We prove that the focal varieties of a $C^{2}$ transnormal function are smooth submanifolds and each regular level set is a tube over either of the focal varieties.

\ni{\bf Key words:}
transnormal function; focal variety; tube; submanifold; Morse index.}

\ni{\it Mathematics Subject Classification (2010):} 53C60, 53C42, 34D23.}
\parskip .001 truein\baselineskip 6pt \lineskip 6pt
\section{Introduction}

\ \ \ \ In Riemannian geometry, the study on isoparametric hypersurfaces has a long history. In the 1930s, Levi-Civita was the first to give the definition of isoparametric functions and isoparametric hypersurfaces\cite{L}. B. Segre gave a complete classification of isoparametric hypersurfaces in an Euclidean space\cite{S}. E. Cartan began to study the isoparametric hypersurfaces in real space forms with constant sectional curvature $c$ systematically\cite{C}. Then many mathematicians started working on isoparametric functions and the classification of isoparametric hypersurfaces in a space form\cite{C1,TP,CS1}. The basic work of isoparametric functions and isoparametric hypersurfaces on general Riemannian manifolds was established by Q.M.Wang who proved
\begin{thmA}~\cite{QM}
Let $M$ be a connected, complete, smooth Riemannian manifold and $f$ a transnormal function on $M$. Then\\
(1)The focal varieties of $f$ are smooth submanifolds of $M$,\\
(2)Each regular level set of $f$ is a tube over either of the focal varieties.
\end{thmA}
Then, Z.Z. Tang and J.Q. Ge made a further study of isoparametric functions and isoparametric hypersurfaces on general Riemannian manifolds\cite{GT,GT1} and gave the relationship between focal sets and focal varieties of transnormal functions.

In Finsler geometry,~\cite{HYS} introduced the concept of isoparametric functions and transnormal functions. Let $(N,F,d\mu)$ be an
$n$-dimensional Finsler manifold with volume form $d\mu$. A $C^{1}$ non-constant function $f$
on $(N,F,d\mu)$ is called \textit{isoparametric} if it is smooth on $N_{f}$ and there are functions $a(t)$ and $b(t)$ such that
\begin{equation}\label{0.1} \left\{\begin{aligned}
&F^{2}(\nabla f)=a(f),\\
&\Delta f=b(f),
\end{aligned}\right.
\end{equation}
where $\nabla f$ denotes the gradient of $f$, which is defined by means of the Legendre transformation, and $\Delta f$ is a nonlinear Finsler-Laplacian of $f$. A function $f$ satisfying only the first equation of (\ref{0.1}) is called transnormal.

Similar to Riemannian geometry, the study of Finsler isoparametric functions and isoparametric hypersurfaces started from Finsler space forms. The flag curvature of a Finsler manifold is a natural generalization of the sectional curvature and we call a (forward) complete and simply connected Finsler manifold with constant flag curvature a (forward) Finsler space form. There are many kinds of special (forward) Finsler space forms, such as Minkowski spaces(with zero flag curvature), Funk spaces(with negative constant flag curvature) and Randers spheres(with positive constant flag curvature).~\cite{HYS,HYS1} classified isoparametric hypersurfaces in Minkowski spaces and Funk spaces.~\cite{HD,HDY} got a complete classification of $d\mu_{BH}-$isoparametric hypersurfaces in a Randers space form and a Funk-type space.~\cite{HCY} studied isoparametric hypersurfaces in Finsler space forms by considering focal points, tubes and parallel hypersurfaces of submanifolds.~\cite{CH} gave the global expression of geodesics, isoparametric functions and isoparametric families on a Randers sphere. In addition to Finsler space forms, there are few results on isoparametric functions and isoparametric hypersurfaces in other cases. M. Xu discussed the correspondences of local isoparametric functions and geodesics by homothetic navigation\cite{XM}. Then he used the same technique to study isoparametric hypersurfaces induced by homothetic navigation in Lorentz Finsler geometry\cite{XM1}.

In order to give a further study on isoparametric functions and isoparametric hypersurfaces in Finsler geometry, in this papar, we study transnormal functions and focal varieties in Finsler manifolds. Unlike the Riemannian case, besides the regularity of isoparametric functions is worse, the Finsler metrics are not necessarily reversible, that is, the geodesics are not necessarily reversible, which cause more difficulties for the research.

For a connected Finsler manifold $(N,F)$, denote the focal varieties of $f$ by $V=\{p\in N\big|\nabla f(p)=0\}$ and $M_{\pm}$ the set where $f$ attains its global maximum value or global minimum value (if either of them exists). As a generalization of Theorem A, set $J=f(N)$, we have

\begin{theo} \label{thm0}
Let $(N,F)$ be a complete, connected Finsler manifold and $f$ a $C^{2}$ transnormal function on $N$. If $a(t)\neq0$ for $t\in\textmd{int}J$ and $a'(t)\neq0$ for $t\in\partial J$, then we have\\
(1)each focal variety of $f$ is a smooth submanifold of $N$.\\
(2)each regular level set of $f$ is a tube over either of the focal varieties with respect to $\overleftarrow{F}$.
\end{theo}

\begin{rema}\label{rema1}
From the condition in Theorem \ref{thm0}, we obviously have $V=M_{\pm}$. Compare with Theorem A, there are more conditions on $f$ and $a(t)$ in Theorem \ref{thm0} which hold in Riemannian case naturally. There exist functions which satisfy Theorem \ref{thm0}(See Example \ref{exam1}). Furthermore, from Lemma \ref{lem3} and Lemma \ref{lem4}, we know that if $a(t)\in C^{2}(\textmd{int}J)$, then $a(t)\neq0$ for $t\in\textmd{intJ}$ and if $a(t)\in C^{2}(J)$, then $a'(t)\neq0$ for $t\in\partial J$. There also exist functions which satisfy $a(t)\in C^{2}(J)$ and $f\in C^{2}(N)$(See Example \ref{exam2}).
\end{rema}

\section{Preliminaries}
~~~~In this section, we will give some definitions and lemmas that will be used in the proof of our main results.
\subsection{Finsler manifolds}
~~~~Let $N$ be an $n$-dimensional manifold and let $TN=\cup_{x\in N}T_xN$ be the tangent bundle of $N$, where $T_xN$ is the tangent space at
$x\in N$. From now on, we use the convention of index ranges $1\leq i,j\leq n$, $1\leq a,b\leq n-1$.

A Finsler metric is a Riemannian metric without quadratic restriction. Precisely, a function $F(x,y)$ on $TN$ is called a \textit{Finsler metric} on a manifold $N$ with local coordinates $(x,y)$, where $x=(x^i)$ and $y=y^i\frac{\partial}{\partial x^{i}}$ ,
if it has the following properties:

(i)\ \ Regularity:\ \ $F(x,y)$ is $C^{\infty}$ on $TN\backslash\{0\}$;

(ii)\ \ Positive homogeneity:\ \ $F(ty)=tF(y),\ \forall t>0, y\in T_xN$;

(iii)\ \ Strong convexity:\ \ The $n\times n$ matrix $(\frac{\partial^2F^2}{\partial y^i \partial y^j}(x,y))(y\neq 0)$ is positive definite.

The fundamental form~$g$ of~$(N,F)$ is given by
\begin{equation*}
g=g_{ij}(x,y)dx^{i} \otimes dx^{j}, ~~~~~~~g_{ij}(x,y)=\frac{1}{2}[F^{2}] _{y^{i}y^{j}}.
\end{equation*}

The projective sphere bundle of~$(N,F)$ is defined by
$$SN=\{y\in TN\setminus\{0\}\big|F(y)=1\}.$$
For a given $x\in N$, the fiber $S_{x}N$ is called the \textit{indicatrix} of $(N,F)$ at $x$.

The projection~$\pi: TN\rightarrow N$ gives rise to the pull-back bundle~$\pi^{\ast}TN$. There exists the unique
\emph{Chern connection}~$\nabla$ with~$\nabla
\frac{\partial}{\partial x^i}=\omega_{i}^{j}\frac{\partial}{\partial
x^{j}}=\Gamma^{i}_{jk}dx^k\otimes\frac{\partial}{\partial x^{j}}$ on the pull-back bundle~$\pi^{\ast}TN$. For $X=X^{i}\frac{\partial}{\partial x^{i}}\in\Gamma(TN)$, the covariant derivative of $X$ along $v=v^{i}\frac{\partial}{\partial x^{i}}\in T_{x}N$ with respect to a reference vector $w\in T_{x}N\setminus\{0\}$ is defined by
$$\nabla^{w}_{v}X(x)=\{v^{j}\frac{\partial X^{i}}{\partial x^{j}}(x)+\Gamma^{i}_{jk}(w)v^{j}X^{k}(x)\}\frac{\partial}{\partial x^{i}}.$$
The equation of geodesics can be expressed by  $\nabla^{\dot{\gamma}}_{\dot{\gamma}}\dot{\gamma}\equiv0$. If a vector field $J(t)$ along $\gamma(t)$ satisfies
\begin{align}\label{B}
\nabla_{\dot{\gamma}}^{\dot{\gamma}}\nabla_{\dot{\gamma}}^{\dot{\gamma}}J+\textbf{R}_{\dot{\gamma}}(J)=0,
\end{align}
then $J(t)$ is called the \textit{Jacobi field} along $\gamma(t)$, where $\textbf{R}_{y}$ is the \textit{flag curvature tensor} which can be expressed by \emph{geodesic coefficients}
\begin{align*}
G^{i}=\frac{1}{4}g^{il}\left\{[F^{2}]_{x^{k}y^{l}}y^{k}-[F^{2}]_{x^{l}}\right\}
\end{align*}
and
$$R^{i}_{k}(y)=2\frac{\partial G^{i}}{\partial x^{k}}-y^{j}\frac{\partial G^{i}}{\partial x^{j}\partial y^{k}}+2G^{j}\frac{\partial G^{i}}{\partial y^{j}\partial y^{k}}-\frac{\partial G^{i}}{\partial y^{j}}\frac{\partial G^{j}}{\partial y^{k}}.$$

Let~${\mathcal L}:TN\rightarrow T^{\ast}N$ denote the \textit{Legendre transformation}, satisfying~${\mathcal L}(\lambda
y)=\lambda {\mathcal L}(y)$ for all~$\lambda>0,~y\in TN$. Moreover, we know from~\cite{SZ} that
\begin{align*}
\mathcal L^{-1}(\xi)=F^{*}(\xi)[F^{*}]_{\xi^{i}}(\xi)\frac{\partial}{\partial x^{i}},\ \ \forall \xi\in T^{*}N\setminus \{0\}, \ \ \mathcal L^{-1}(0)=0,
\end{align*}
where $F^{*}$ is the dual metric of $F$. For a smooth function~$f: N\rightarrow \mathbb{R}$, the \textit{gradient vector} of~$f$ at~$x$ is defined as~$\nabla f(x)={\mathcal
L}^{-1}(df(x))\in T_{x}N$. Hessian of $f$ is defined by
\begin{align}\label{2.3}
\textrm{Hes}f(X,Y):=g_{\nabla f}(\nabla^{2}f(X),Y)=X(Yf)-(\nabla^{\nabla f}_{X}Y)f,
\end{align}
where $X,Y\in T_{x}N$, $\nabla^{2}f(x)=\nabla^{\nabla f}_{X}\nabla f\in T_{x}^{*}N\otimes T_{x}N$ for~$x\in N_{f}=\{x\in N|df(x)\neq 0\}$. From \cite{YHS}, $\textrm{Hes} f(X,Y)=\textrm{Hes} f(Y,X)$.
\subsection{Anisotropic submanifolds of a Finsler manifold}
~~~~Let $(N,F)$ be an $n$-dimensional Finsler manifold and $\phi: M\to(N, F)$ be an $m$-dimensional immersion. For simplicity, we will denote $d\phi(x)$ by $x$ and $d\phi X$ by $X$. Let
$$\mathcal{V}(M)=\{(x,\xi)~|~x\in M,\xi\in T_x^{*}N,\xi (X)=0,\forall X\in T_xM\},$$
which is called the
\textit{normal bundle} of $M$ \cite{SZ1}. Set $\mathcal{N}M={\mathcal L}^{-1}(\mathcal{V}(M))\subset TN$. In general, $\mathcal{N}M$ is not a vector bundle. For any locally defined \textit{normal vector} $\textbf{n}\in\mathcal{N}_{x}M$ which satisfies $g_{\textbf{n}}(\textbf{n},X)=0$ for any $X\in T_{x}M$, we can define a Riemannian metric $\hat{g}:=\phi^{*}g_{\textbf{n}}$. We call $(M,\hat{g})$ an
\textit{anisotropic submanifold} of $(N,F)$. For any $X\in T_xM$, the \emph{shape operator}~${A}_{\mathbf{\textbf{n}}}:T_xM\rightarrow T_xM$ is defined by
\begin{equation}\label{A}
{A}_{\mathbf{\textbf{n}}}X=-\left(\nabla^{\mathbf{\textbf{n}}}_{X}\mathbf{\textbf{n}}\right)^{\top}_{g_{\mathbf{\textbf{n}}}}.
\end{equation}
Furthermore, we denote the unit normal bundle by
$$\mathcal{V}^0(M)=\{\nu\in \mathcal{V}(M)|F^*(\nu)=1 \}.$$
Set $\mathcal{N}^0M={\mathcal L}^{-1}(\mathcal{V}^0(M))$, $\textbf{n}\in \mathcal{N}^{0}_{x}M$ denotes the \textit{unit normal vector} of $M$.

If $m=n-1$, there exists a global normal vector field $\textbf{n}$. Then $(M,\hat{g})$ is an oriented Riemannian manifold which is called an oriented $\textit{anisotropic hypersurface}$. In this case,
\begin{equation*}
{A}_{\mathbf{\textbf{n}}}X=-\nabla^{\mathbf{\textbf{n}}}_{X}\mathbf{\textbf{n}}.
\end{equation*}
\subsection{Focal sets and tubes}
~~~~Define the normal exponential map $E:\mathcal{N}M\rightarrow N$ satisfying $E(x,\textbf{n})=\exp_{x}\textbf{n}$, which is $C^{\infty}$ on $TN\setminus\{0\}$. The \textit{focal points} of $M$ are the critical values of the normal exponential map $E$. The focal point $p\in N$ has multiplicity $m$ means that $E_{*}$ at $(x,\textbf{n})$ has nullity $m>0$.

For the immersion $\phi:M\rightarrow (\overrightarrow{N},F)$ with codimension $n-m\geq1$, if $n-m>1$, we define $\phi_{s}M(s>0)$ by the map $\phi_{s}:\mathcal{N}^{0}M\rightarrow N$ satisfying $\phi_{s}(x,\textbf{n})=E(x,s\textbf{n})$. If $(x,s\textbf{n})$ is not a critical point of $E$, then $\phi_{s}$ is an $(n-1)$-dimensional immersion in $(N,F)$. From \cite{HCY}, in the neighborhood $U$ of $x\in\phi(M)$, there exists a local tube over $U$ belonging to $\phi_{s}M$ for all sufficiently small $s>0$. We call $\phi_{s}M$ a \textit{tube} over $M$.

If $M$ is a hypersurface, then $\mathcal{N}^{0}M$ is a double covering of $M$. In this case, for a given unit normal field $\textbf{n}$ of $M$, define $\phi_{s}:M\rightarrow (\overrightarrow{N},F)$ by
\begin{align}\label{0}
\phi_{s}(x)=E(x,s\textbf{n}(x)).
\end{align}
$\phi_{s}$ is a local immersion for sufficiently small $s$. If there is no focal point of $M$ on $\phi_{s}M$, $\phi_{s}M$ is a \textit{parallel hypersurface} over $M$ as $s>0$ and $M$ is a \textit{parallel hypersurface} over $\phi_{s}M$ as $s<0$ with respect to $\textbf{n}$. For simplicity, we call that $\phi_{s}M$ is a \textit{parallel hypersurface} over $M$.

\subsection{The reverse metric}
~~~~Let $(N, F)$ be an oriented smooth Finsler manifold. The reverse metric of $F$ is defined by
$
\overleftarrow{F}(x,y)=F(x,-y)$, where $y\in T_{x}N$, $x\in N$. Then $
\overleftarrow{F}^*(\xi)=F^*(-\xi)
$, where $\xi\in T^*_{x}N$.
Obviously, by direct computation,
\begin{align*}
\overleftarrow{g}_{ij}(y)=g_{ij}(-y),~~\overleftarrow{\Gamma}^{i}_{jk}(y)
=\Gamma^{i}_{jk}(-y),~~\overleftarrow{N}^{i}_{j}(y)=-N^{i}_{j}(-y).
\end{align*}
Let $f: N \to \mathbb{R}$ be a $C^{1}$ non-constant function on  $(N, F)$. Then from the definition of gradient vector of $f$, we have $$\overleftarrow{\nabla}(-f)=-\nabla f.$$
Then we know if $f$ is a transnormal function with respect to $F$ if and only if $-f$ is a transnormal function with respect to $\overleftarrow{F}$ which satisfies
\begin{equation*}
\overleftarrow{F}(\overleftarrow{\nabla}(-f))= a(f).
\end{equation*}
Furthermore, for the level set $M_{t}$, $\textbf{n}=\frac{\nabla f}{F(\nabla f)}$ and
$\overleftarrow{\textbf{n}}=\frac{\overleftarrow{\nabla} (-f)}{\overleftarrow{F}(\overleftarrow{\nabla}(-f))}$ are the unit normal vector fields of $M$.

\section{Some properties of transnormal functions}
~~~~In this section, let $f$ be a transnormal function on a connected and forward complete Finsler manifold $(N,F)$, which is $C^{1}$ on $N$ and smooth on $N_{f}$. Recall that a geodesic segment $\gamma(s)$ is called an \textit{$f-$segment} if $f(\gamma(s))$ is an increasing function of $s$, where $s$ is an arc-length parameter.

\begin{lemm}~\cite{HYS} \label{lem1}
We have the following results on $\textmd{int}J$.\\
(1)There is a function $\rho$ defined on $\textmd{int}J$ such that $\nabla\rho=\frac{\nabla f}{F(\nabla f)}$.\\
(2)The integral curves of the gradient vector field $\nabla f$ in $\textmd{int}J$ are all geodesics($f$-segments).\\
(3)The regular level hypersurfaces $M_{t}=f^{-1}(t)$ are parallel along the direction of $\nabla f$.\\
(4)$d(x,M_{t_{2}})=d(M_{t_{1}},y)=\int^{t_{2}}_{t_{1}}\frac{1}{\sqrt{a(t)}}dt$ for any $[t_{1},t_{2}]\in \textmd{int}J$, $x\in M_{t_{1}}$, $y\in M_{t_{2}}$ and the $f$-segments are the shortest curves among all curves connecting $M_{t_{1}}$ and $M_{t_{2}}$.
\end{lemm}

\begin{lemm}\label{lem2}
Suppose that $\beta$ is the only critical value of $f$ in $[\alpha, \beta]\subseteq J$, then
\begin{equation}\label{3.2}
d(M_{\alpha}, M_{\beta})=\int^{\beta}_{\alpha}\frac{df}{\sqrt{a(f)}}
\end{equation}
and the improper integral in (\ref{3.2}) converges.
\end{lemm}
\proof
Let $\{p\}$ be any point in $M_{\beta}$. Since $\{f\geqslant \beta\}$ in $N$ is closed and $N$ is connected, then there are points $\{p_{i}\}\in \{f<\beta\}$ in $N$ such that $\lim\limits_{i\rightarrow\infty}d(p,p_{i})=0$. Set $f(p_{i})=c_{i}$.

Centered at $p$, there exists a normal geodesic $\gamma(s)=\exp_{p} sX$ through $p_{i}$ which satisfies $\gamma(0)=p$, $\gamma'(0)=X$, where $X\in S_{p}N$. Let $U(p)$ be a compact neighborhood covering the geodesic sphere determined by $\gamma(s)$ with radius less than $\max d(p,p_{i})$ and $\lambda_{0}=\max\limits_{U(p)}\frac{\overleftarrow{F}(y)}{F(y)}$. Then we have
\begin{align}\label{3.3}
d(p_{i},p)\leq L(\overleftarrow{\gamma}(s))=\int^{s_{i}}_{0}|F(-\gamma'_{i}(s))|ds=\int^{s_{i}}_{0}|\overleftarrow{F}(\gamma'_{i}(s))|ds \leq \lambda_{0}\int^{s_{i}}_{0}|F(\gamma'_{i}(s))|ds.
\end{align}
Hence,
\begin{align}\label{3.4}
d(p_{i},p)\leq \lambda_{0}s_{i}=\lambda_{0}d(p,p_{i}).
\end{align}
Namely,
\begin{align}\label{3.5}
\lim_{i\rightarrow\infty}d(p_{i},p)=0.
\end{align}

Let $\sigma_{i}(t)$ be any piecewise $C^{1}$ curve which goes from $M_{\alpha}$ to $p$ through $p_{i}$. For every $p_{i}\in(\alpha,\beta)$, $\sigma_{1_{i}}$ is a curve from $M_{\alpha}$ to $M_{p_{i}}$ and $\sigma_{2_{i}}$ is a curve from $M_{p_{i}}$ to $p$ such that $\sigma_{i}=\sigma_{1_{i}}\cup\sigma_{2_{i}}$. Furthermore, $\inf L(\sigma_{i})=\inf L(\sigma_{1_{i}})+\inf L(\sigma_{2_{i}})$. Combine (\ref{3.5}), we obtain
\begin{align*}
d(M_{\alpha},M_{\beta})=\lim_{i\rightarrow\infty}d(M_{\alpha},M_{p_{i}}).
\end{align*}
This completes the proof.
\endproof

\begin{lemm}\label{lem3}
$f$ has no critical value(i.e. $a(t)\neq 0$) in $\textmd{int}J$ when $a(t)\in C^{2}(\textmd{int}J)$.
\end{lemm}
\proof
Suppose $\beta$ the only critical value of $f$ in $\textmd{int}J$, then we have $a(\beta)=0$ and $a'(\beta)=0$. If $t\in[\beta-\varepsilon,\beta]\subseteq\textmd{int}J$, from Taylor's theorem,
\begin{align}\label{101}
a(t)=a(\beta)+a'(\beta)(t-\beta)+\frac{1}{2}a''(\xi)(t-\beta)^2=\frac{1}{2}a''(\xi)(t-\beta)^2,\ \ \ \ \xi\in(t,\beta).
\end{align}
Since $a(t)\in C^{2}([\beta-\varepsilon,\beta])$, there exists a constant $K$ such that $|a''(\xi)|\leq 2K$. From (\ref{0.1}) and (\ref{101}), $0\leq a(t)\leq K(t-\beta)^2$, then the integral in (\ref{3.2}) would diverge. There exists a contradiction with Lemma \ref{lem2}.
\endproof

Combine Lemma \ref{lem2} and $a(t)\neq 0$ on $\textmd{int}J$, $s(t)=\int_{t_{0}}^{t}\frac{1}{\sqrt{a(x)}}dx$ is convergent, where $t_{0},t\in J$. Hence, from (\ref{0}), set $r$ the distance from $M_{c}$ to $M_{\pm}$, the smooth map $\phi_{r}:M_{c}\rightarrow N$ is
\begin{align}\label{4.3}
\phi_{r}(x)=\exp_{x}(r\textbf{n}).
\end{align}
It's easy to show that the focal points of $M_{c}$ are also the critical values of $\phi_{r}$. Hence, $\phi_{r}$ is called the \textit{focal map} of $M_{c}$. Suppose $\max J=\beta<\infty$ if the maximal value of $J$ exists, $M_{+}=f^{-1}(\beta)=\phi_{r}M_{c}$. Note that fix a $c$ in $\textmd{int}J$, the hypersurface $M_{c}$ may not be connected.
Define $\textmd{Foc}_{t}(f)=\{\phi_{s(t)}(x)\in N\big|d\phi_{s(t)}X=0, X\in T_{x}M_{c}\setminus\{0\}\}$ the focal set of the level set $M_{t}=\phi_{s(t)}M=f^{-1}(t)$. Since $\phi_{s(t)}:M\rightarrow M_{t}$ is a diffeomorphism if $s<r$, we easily get $\textmd{Foc}(f)\subseteq M_{\pm}$. For two different level hypersurfaces $M_{t_{1}}$ and $M_{t_{2}}$, from Lemma \ref{lem2}, we have $\textmd{Foc}_{t_{1}}(f)=\textmd{Foc}_{t_{2}}(f)$. Hence, we can simply write $\textmd{Foc}_{t}(f)$ as $\textmd{Foc}(f)$.
\begin{lemm}\label{lem4}
If $a(t)\in C^{2}(J)$, then $a'\neq0$ on $\partial J$. Moreover, $a'<0$(resp.$a'>0$) at the maximal(resp.minimal) point of $J$ if it exists.
\end{lemm}
\proof
Suppose $\beta$ the only critical value of $f$ in $[\beta-\varepsilon,\beta]$ for small $\varepsilon$. From Taylor's theorem,
\begin{align*}
a(t)=a(\beta)+a'(\beta)(t-\beta)+\frac{1}{2}a''(\xi)(t-\beta)^2=a'(\beta)(t-\beta)+\frac{1}{2}a''(\xi)(t-\beta)^2,
\end{align*}
where $t\in(\beta-\varepsilon,\beta)$. Firstly, we can prove $a'(\beta)\neq0$. If $a'(\beta)=0$, similar to the analysis in Lemma \ref{lem3}, there exists a contradiction with Lemma \ref{lem2}. Moreover, we claim $a'(\beta)<0$. If $a'(\beta)>0$, we have $a(t)<0$ due to the continuity of $a(t)$. This is in contradiction with $a(t)\geq0$.
\endproof

From Lemma \ref{lem1}, Lemma \ref{lem2} and the first variation formula for arc length, we have
\begin{lemm}\label{lem5}
(1) Given a point $x$ in $f^{-1}(\textmd{int}J)$, there is a unique maximal $f$-segment through $x$.\\
(2) Every maximal $f$-segment intersects each regular level set of $f$ exactly once and orthogonally.\\
(3) Let $\gamma(s)$ be an $f$-segment, then for $s_{1}<s_{2}$, $d(M_{f(\gamma(s_{1}))},M_{f(\gamma(s_{2}))})=s_{2}-s_{1}$.\\
(4) The $f$-segment $\gamma(s)$ is also perpendicular to $M_{+}$(resp.$M_{-}$).\\
\end{lemm}
\proof
(1),~(2),~(3) are direct corollaries of Lemma \ref{lem1}, Lemma \ref{lem2}. Here, we focus on the proof of (4).

Set $\sigma(t)$ any curve on $M_{+}$ through $p$ satisfying $\sigma(0)=p$. $\gamma(s)$ is a geodesic which intersects $M_{c}$ and $M_{+}$ at $x=\gamma(0)$ and $p=\gamma(r)$, respectively. Let $\widetilde{\sigma}(t)$ be any curve on $M_{c}$ through $x$ satisfying $\widetilde{\sigma}(0)=x$. Obviously, $\langle \widetilde{\sigma}'(0),\gamma'(0)\rangle=0$. We only need to prove $\langle \sigma'(0), \gamma'(r)\rangle$=0.

Set $\Phi(s,t)$ any smooth variation which satisfies $\Phi(s,0)=\gamma(s)$, $\Phi(0,t)=\widetilde{\sigma}(t)$ and $\Phi(r,t)=\sigma(t)$. Set $V(s)=\frac{\partial\Phi}{\partial t}\big|_{t=0}$ and $T(s)=\frac{\partial\Phi}{\partial s}\big|_{t=0}$ which satisfy
$$V(0)=\widetilde{\sigma}'(0),\ \ \ \ V(r)=\sigma'(0),\ \ \ \ T(0)=\gamma'(0),\ \ \ \ T(r)=\gamma'(r),\ \ \ \ \langle T(0), V(0)\rangle=0.$$
From the first variation formula for arc length, we know
\begin{align*}
0=L'(0)=\langle T(s), V(s)\rangle\big|^{s=r}_{s=0}-\int^{r}_{0}g_{T(s)}(\nabla^{T(s)}_{T(s)}T(s), V(s))ds=\langle V(r), T(r)\rangle=\langle \sigma'(0), \gamma'(r)\rangle.
\end{align*}
This completes the proof of (4).
\endproof
\section{Focal varieties and tubes}
\subsection{Hessian}
~~~~From now on, suppose $f\in C^{2}(N)$. Hessian on $N_{f}$ is defined by (\ref{2.3}). Now, we need to give the definition of Hessian on focal varieties which satisfies $\nabla f=0$. In this case, $g_{\nabla f}$ is not defined. For $p\in N$, $U_{p}(x^{i})$ is a local coordinate system at $p$. For the convenience of writing, we use $f_{ij}$ and $f_{k}$ instead of $\frac{\partial^{2}f}{\partial x^{i}\partial x^{j}}$ and $\frac{\partial f}{\partial x^{k}}$, respectively.

For any $p\in M_{+}$, $X_{p},Y_{p}\in T_{p}M_{+}$, extend $X_{p},Y_{p}$ to $U_{p}$. Set $q\in U_{p}\setminus M_{+}$. Define
\begin{align*}
\textmd{Hes}f(X_{p},Y_{p})&=\lim_{q\rightarrow p}\textmd{Hes}f(X_{q},Y_{q})\\
&=\lim_{q\rightarrow p}X^{i}_{q}Y^{j}_{q}f_{ij}(q)-f_{k}(q)\Gamma^{k}_{ij}(\nabla f)X^{i}_{q}Y^{j}_{q}.
\end{align*}
Since $\lim\limits_{q\rightarrow p}f_{k}(q)=0$ and $\Gamma^{k}_{ij}(\nabla f)$ is homogenous of degree $0$ which is bounded in $U_{p}\setminus M_{+}$, we have
\begin{align}\label{4.1}
\textmd{Hes}f(X_{p},Y_{p})=\lim_{q\rightarrow p}X^{i}_{q}Y^{j}_{q}f_{ij}(q)=X^{i}_{p}Y^{j}_{p}f_{ij}(p).
\end{align}
Obviously, $\textmd{Hes}f(X_{p},Y_{p})$ is self-adjoint. For the convenience of writing, we will use $H$ instead of $\textmd{Hes}$ in the following. Then, for $\xi=\mathcal{L}(X)$, we define
\begin{align}\label{4.4}
H^{*}(\xi)(Y)=H(X,Y).
\end{align}
Using $H$ and $H^{*}$, for any $p\in N$, denote
$$E_{0}(p)=\{X\in T_{p}N\big|~H(X)=0\},$$
$$E^{*}_{0}(p)=\{\xi\in T^{*}_{p}N\big|~H^{*}(\xi)=0\,\},$$
$$E_{\mu}^{*}(p)=\{\xi\in T^{*}_{p}N\big|~H^{*}(\xi)=\mu\xi\},$$
$$E^{\bot}_{0}(p)=\{\xi\in T^{*}_{p}N\big|~\xi(X)=0,X\in E_{0}(p)\}$$
and $\mathcal{U}E$ as the unit vectors in a vector space $E$ with respect to $F$(or $F^{*}$).

\begin{remark}
Since $H^{*}$ is not necessarily a linear operator, $E^{*}_{\mu}$ is not necessarily a vector space.
\end{remark}

Set $\gamma_{x}(s)=\exp_{x}s\textbf{n}$ a normal geodesic in $N$, where $x=\gamma_{x}(0)\in M_{c}$ and $\textbf{n}=\gamma'_{x}(0)$ is the unit normal vector of $M_{c}$. Define the map $\sigma^{*}:M_{c}\rightarrow T^{*}N$ which satisfies $\sigma^{*}(x)=\mathcal{L}(\gamma'_{x}(r))$. For $p\in M_{+}$, set $A_{p}=\phi^{-1}_{r}(p)\subseteq M_{c}$. The relationship between $\sigma^{*}A_{p}$, $\mathcal{U}E^{*}_{\lambda}(p)$ and $E^{\bot}_{0}(p)$ can be obtained.
\begin{lemm}\label{lemX}
(1)~$\sigma^{*}A_{p}\subseteq \mathcal{U}E^{*}_{\lambda}(p)$, where $\lambda=\frac{1}{2}a'(\beta)$.\ \ (2)~$E^{*}_{\lambda}(p)\subseteq E^{\bot}_{0}(p)$.
\end{lemm}
\proof
(1) For any $\tau\in \sigma^{*}A_{p}$, there exist $x\in A_{p}$ and $\textbf{n}\in \mathcal{N}_{x}^{0}M_{c}$ such that $\gamma_{x}(s)=\exp_{x}s\textbf{n}$ satisfying $p=\gamma_{x}(r)$ and $\eta=\mathcal{L}^{-1}(\tau)=\gamma_{x}'(r)$. Let $\tilde{U}(p)$ be a neighborhood of $p$ in $N$ and $T(s)=\gamma_{x}'(s)$. Extend $T$ to the vector field $\widetilde{T}$ in $\widetilde{U}(p)\subseteq N$ such that $\widetilde{T}=\frac{\nabla f}{\sqrt{a}}$ in $\widetilde{U}(p)\setminus M_{+}$ and $\widetilde{T}\big|_{\gamma(s)}=T$. For any $Y\in TN$, $q\in\widetilde{U}(p)\setminus M_{+}$, from (\ref{4.4}),
\begin{align*}
H^{*}(\tau)(Y)&=H(\eta,Y)=\lim\limits_{q\rightarrow p}g_{\widetilde{T}}(\nabla_{\widetilde{T}}^{\widetilde{T}}(\sqrt{a}\widetilde{T}),Y)\\
&=\lim\limits_{q\rightarrow p}\widetilde{T}(\sqrt{a})g_{\widetilde{T}}(\widetilde{T},Y)+\lim\limits_{q\rightarrow p}\sqrt{a}g_{\widetilde{T}}(\nabla^{\widetilde{T}}_{\widetilde{T}}\widetilde{T},Y)\\
&=\lim\limits_{q\rightarrow p}\frac{1}{2\sqrt{a}}a'g_{\widetilde{T}}(\widetilde{T},\nabla f)g_{\widetilde{T}}(\widetilde{T},Y)=\frac{1}{2}a'(\beta)\tau(Y).
\end{align*}
Hence, $H^{*}(\tau)=\lambda\tau$, where $\lambda=\frac{1}{2}a'(\beta)$.

(2) For any $\xi\in E^{*}_{\lambda}(p)$, we have $H^{*}(\xi)=\lambda\xi$, where $\lambda\neq0$. For any $X\in E_{0}$, we have
\begin{align*}
0=H^{*}(\xi)(X)=H(\mathcal{L}^{-1}(\xi),X)=\lambda\xi(X).
\end{align*}
Hence, $\xi(X)=0$. Then we have $\xi\in E^{\bot}_{0}(p)$. Namely, $E^{*}_{\lambda}(p)\subseteq E^{\bot}_{0}(p)$.
\endproof
\subsection{Tubes of smooth submanifolds}
~~~~Let $U$ be an open subset in $M_{c}$ such that rank$(\phi_{r}\big|_{U})=k$ and $W=\phi_{r}(U)$ a smooth $k$-dimentional submanifold of $N$. Set $x\in U$, $p=\phi_{r}(x)$. In this section, we denote $A_{p}=\phi^{-1}_{r}(p)\cap U$. From Lemma \ref{lem5}, we obviously know that $\sigma^{*}$ maps $A_{p}$ to $\mathcal{V}^{0}_{p}W$. Consider the relationship between $T_{p}W$, $E_{0}(p)$, $\mathcal{V}^{0}_{p}W$ and $\mathcal{U}E^{\bot}_{0}(p)$, we have
\begin{lemm}\label{lem7}
(1)~$T_{p}W\subseteq E_{0}(p)$,\ \ (2)~$\mathcal{V}^{0}_{p}W=\mathcal{U}E^{\bot}_{0}(p)$.
\end{lemm}
\proof
(1)~Set $\forall X\in T_{p}W$. Extend $X$ to $\widetilde{U}(p)$ such that $g_{\widetilde{T}}(\widetilde{T},X)=0$, where $\widetilde{T}$ and $\widetilde{U}(p)$ are defined in the proof of (1) in Lemma \ref{lemX}. For any $Y\in TN$ and $q\in \widetilde{U}(p)\setminus M_{+}$, from (\ref{4.1}),
\begin{align*}
H(X,Y)(p)&=\lim_{q\rightarrow p}H(X,Y)(q)=\lim_{q\rightarrow p}g_{\widetilde{T}}(\nabla^{\widetilde{T}}_{X}(\sqrt{a}\widetilde{T}),Y)\\
&=\lim_{q\rightarrow p}g_{\widetilde{T}}(\nabla^{\widetilde{T}}_{Y}(\sqrt{a}\widetilde{T}),X)\\
&=\lim_{q\rightarrow p}Y(\sqrt{a})g_{\widetilde{T}}(\widetilde{T},X)+\lim_{q\rightarrow p}\sqrt{a}g_{\widetilde{T}}(\nabla^{\widetilde{T}}_{Y}\widetilde{T},X)\\
&=\lim_{q\rightarrow p}\sqrt{a}g_{\widetilde{T}}(\nabla^{\widetilde{T}}_{Y}\widetilde{T},X)
\end{align*}
Since
\begin{align}\label{4.2}
g_{\widetilde{T}}(\nabla^{\widetilde{T}}_{Y}\widetilde{T},X)+g_{\widetilde{T}}(\widetilde{T},\nabla^{\widetilde{T}}_{Y}X)+2C_{\widetilde{T}}(\nabla^{\widetilde{T}}_{Y}\widetilde{T},\widetilde{T},X)=Yg_{\widetilde{T}}(\widetilde{T},X),
\end{align}
we have
$$\lim_{q\rightarrow p}g_{\widetilde{T}}(\nabla^{\widetilde{T}}_{Y}\widetilde{T},X)=-\lim_{q\rightarrow p}g_{\widetilde{T}}(\widetilde{T},\nabla_{Y}^{\widetilde{T}}X)=-g_{\eta}(\eta,\nabla_{Y}^{\eta}X).$$
In addition, $\lim\limits_{q\rightarrow p}\sqrt{a}=0$. Hence, we have $H(X,Y)(p)=0$ for any $Y\in TN$. Namely, $T_{p}W\subseteq E_{0}(p)$.

(2) From Lemma \ref{lem5}, $\sigma^{*}(A_{p})\subseteq\mathcal{V}^{0}_{p}W$. Combine Lemma \ref{lemX}, $\sigma^{*}(A_{p})\subseteq \mathcal{U}E^{\bot}_{0}(p)$.
Namely, $\sigma^{*}(A_{p})\subseteq\mathcal{V}^{0}_{p}W\cap \mathcal{U}E^{\bot}_{0}(p)$. Since $\sigma^{*}$ is injective and $A_{p}=\phi^{-1}_{r}(p)\cap U$ is $(n-1-k)$-dimensional, $\dim \sigma^{*}(A_{p})=n-1-k$. In addition, $\mathcal{V}_{p}W$ and $E^{\bot}_{0}(p)$ are two $(n-k)$-dimensional linear spaces. Hence, $\mathcal{V}^{0}_{p}W=\mathcal{U}E^{\bot}_{0}(p)$.
\endproof

\begin{rema}
(2) in Lemma \ref{lem7} means $E_{0}(p)=T_{p}W$.
\end{rema}

We can also express $\sigma^{*}(x)=\mathcal{L}(\gamma'_{x}(r))$ as $\sigma^{*}=\mathcal{L}\circ \Phi_{r}\circ \textbf{n}$, where $\textbf{n}(x)$ embeds $M_{c}$ into the unit tangent bundle $SN$ and $x\in M_{c}$. Define $\Phi_{s}:SN\rightarrow SN$ by $\Phi_{s}(x,\textbf{n})=(\exp_{x}s\textbf{n},\frac{d}{ds}\exp_{x}s\textbf{n})$, then we obtain
\begin{lemm}\label{lem6}
$\Phi_{r}:SN\rightarrow SN$ is a diffeomorphism induced by the geodesic flow of $N$ for the fixed length $r$.
\end{lemm}
\proof
Since $\Phi_{r}(x,\gamma'_{x}(0))=(\gamma_{x}(r),\gamma'_{x}(r))$, it is obvious that $\Phi_{r}$ is injective and smooth. For $(\widetilde{x},\widetilde{\textbf{n}})\in SN$, set $\overleftarrow{\gamma}_{\widetilde{x}}(s)=\exp_{\widetilde{x}}s(-\widetilde{\textbf{n}})$, where $\overleftarrow{\gamma}_{\widetilde{x}}(0)=\widetilde{x}$, $\overleftarrow{\gamma}'_{\widetilde{x}}(0)=-\widetilde{\textbf{n}}$. Then for $x=\overleftarrow{\gamma}_{\widetilde{x}}(r)$ and $\textbf{n}=-\frac{d}{ds}\overleftarrow{\gamma}_{\widetilde{x}}(s)\big|_{s=r}$, we have $\Phi_{r}(x,\textbf{n})=(\tilde{x},\tilde{\textbf{n}})$. Therefore, $\Phi_{r}$ is a surjection. This completes the proof.
\endproof

From Lemma \ref{lemX}, Lemma \ref{lem7} and Lemma \ref{lem6},
\begin{lemm}\label{lem8}
$\phi_{r}^{-1}W$ is a tube of radius $r$ over $W$ with respect to $\overleftarrow{F}$.
\end{lemm}
\proof
Since $\mathcal{N}^{0}(M_{c})$ is closed in $SN$, $\Phi_{r}$ and $\mathcal{L}:SN\rightarrow S^{*}N$ are both diffeomorphism, $\sigma^{*}M_{c}$ is closed in $S^{*}N$. Moreover, $\sigma^{*}M_{c}\cap\mathcal{V}^{0}_{p}W$ is closed in $\mathcal{V}^{0}_{p}W$. For any $x\in M_{c}$, set $U_{x}$ an open neighborhood in $N$. Meanwhile, there exists an open neighborhood $V_{x}$ of $\textbf{n}(U_{x})$ in $SN$ such that
$$\sigma^{*}(U_{x}\cap M_{c})\cap\mathcal{V}^{0}_{p}W=\mathcal{L}\circ\Phi_{r}(V_{x})\cap\mathcal{V}^{0}_{p}W.$$
Since $\mathcal{L}\circ\Phi_{r}(V_{x})$ is open in $S^{*}N$, $\sigma^{*}(U_{x}\cap M_{c})\cap\mathcal{V}^{0}_{p}W$ is open in $\mathcal{V}^{0}_{p}W$. Hence, $\sigma^{*}M_{c}\cap\mathcal{V}^{0}_{p}W$ is open in $\mathcal{V}^{0}_{p}W$. From above analysis, $\mathcal{V}^{0}W=\sigma^{*}\phi_{r}^{-1}W$.

For any $(p,\eta)\in\mathcal{N}^{0}W$, set $\overleftarrow{\phi}_{s}(p,-\eta)=\overleftarrow{E}(p,-s\eta)=\overleftarrow{\exp}_{p}s(-\eta)$. In addition, $x=\overleftarrow{\exp}_{p}r(-\eta)\in\phi_{r}^{-1}W\subseteq M_{c}$ is not a critical point of $\overleftarrow{E}$. Hence, $\phi_{r}^{-1}W=\overleftarrow{\phi}_{r}\overleftarrow{\mathcal{N}}^{0}W$ is the tube of radius $r$ over $W$ with respect to $\overleftarrow{F}$.
\endproof

Combine Lemma \ref{lemX} and Lemma \ref{lem8}, we have
\begin{coro}\label{coro1}
$E^{\bot}_{0}(p)=E^{*}_{\lambda}(p)$, where $p\in W$.
\end{coro}
\proof On the one hand, for any $\tau\in \mathcal{U}E^{\bot}_{0}(p)=\mathcal{V}^{0}_{p}W$, from Lemma \ref{lem8}, we know $\tau\in\sigma^{*}\phi_{r}^{-1}p$. Similar to the proof of (1) in Lemma \ref{lemX}, $\sigma^{*}\phi_{r}^{-1}p\subseteq \mathcal{U}E^{*}_{\lambda}(p)$. Namely, $\mathcal{U}E^{\bot}_{0}(p)\subseteq \mathcal{U}E^{*}_{\lambda}(p)$. On the other hand, for any $\tau\in \mathcal{U}E_{\lambda}^{*}(p)$, from (2) in Lemma \ref{lemX}, $\mathcal{U}E^{*}_{\lambda}(p)\subseteq \mathcal{U}E^{\bot}_{0}(p)$. Hence, $\mathcal{U}E^{\bot}_{0}(p)=\mathcal{U}E^{*}_{\lambda}(p)$. This completes the proof.
\endproof

From Lemma \ref{lem8}, we have
\begin{lemm}\label{lem9}
For any $p\in W$, there is a neighborhood $W_{1}$ of~$p$ in $W$ such that $U_{\varepsilon}(W_{1})\cap M_{+}=W_{1}$ for a sufficiently small $\varepsilon>0$, where $U_{\varepsilon}(W_{1})$ is the tubular neighborhood of $W_{1}$ of radius $\varepsilon$ with respect to $F$.
\end{lemm}
\proof
For any $\eta\in\mathcal{N}^{0}_{p}W$, $\tau=\mathcal{L}(\eta)\in \mathcal{V}^{0}_{p}W$, there exists a unique geodesic $\gamma_{\eta}(s)$ determined by $p=\gamma_{\eta}(0)$ and $\eta=\gamma_{\eta}'(0)$ for $F$. Denote the function $f_{\eta}(s)=f(\gamma_{\eta}(s))$, we have $f_{\eta}(0)=\beta$, $f'_{\eta}(0)=df|_{p}(\eta)=0$. Combine Lemma \ref{lem7} and Corollary \ref{coro1},  $\mathcal{U}E^{\bot}_{0}(W)=\mathcal{U}E^{*}_{\lambda}(W)=\mathcal{V}^{0}W$, we have
\begin{align*}
f''_{\eta}(0)=H(\eta,\eta)=H^{*}(\tau)(\eta)=\lambda\tau(\eta)=\lambda.
\end{align*}
From Lemma \ref{lem4} and (1) of Lemma \ref{lemX}, $\lambda<0$. Hence, $f''_{\eta}(0)<0$. Lifting $f$ to $SN$, $f''_{\eta}(s)$ is continuous in $\eta$ and $s$. Hence, there exists a neighborhood $W_{1}$ of $p$ in $W$ and a sufficiently small $\varepsilon$ such that $f''_{\eta}(s)<\frac{\lambda}{2}<0$ for $s\in(0,\varepsilon)$ and $\eta=\mathcal{L}^{-1}(\tau)$, where $\tau\in \mathcal{U}E^{\bot}_{0}(W_{1})$. Hence,
\begin{align*}
f_{\eta}(s)=f_{\eta}(0)+f'_{\eta}(0)s+\frac{1}{2}f''_{\eta}(\xi)s^2<\beta,\ \ \ \ s\in(0,\varepsilon).
\end{align*}
Namely, $f$ is strictly less than $\beta$ if restricted to $U_{\varepsilon}(W_{1})\setminus W_{1}$.
\endproof
\subsection{Maximal rank problem}
~~~~Denote $m=\min\limits_{p\in M_{+}}\{\dim E_{0}(p)\}$, $M_{+}^{m}=\{p\in M_{+}\big|\dim E_{0}(p)=m\}=\{p\in M_{+}\big|\max \textmd{rank} H_{p}=n-m\}$,
$M^{m}_{c}=\phi^{-1}_{r}(M^{m}_{+})$, $\mathring{M}^{m}_{c}$ the subset of $M^{m}_{c}$ where $\phi_{r}$ has maximal rank $m$ and $\mathring{M}^{m}_{+}=\phi_{r}\mathring{M}^{m}_{c}$, we have
\begin{lemm}\label{lem10}
$\mathring{M}^{m}_{c}$ is open and dense in $M^{m}_{c}$.
\end{lemm}
\proof
For any $x\in\mathring{M}^{m}_{c}$, there exists $U(x)\subseteq \mathring{M}^{m}_{c}\subseteq M^{m}_{c}$ in which rank$\phi_{r}\big|_{U(x)}\geq m$. Consider the condition $\phi_{r}$ has maximal rank $m$, we have $\phi_{r}$ has constant rank $m$ on $U(x)$. Hence, $\mathring{M}^{m}_{c}$ is open in $M^{m}_{c}$.

Suppose $\widetilde{U}=\textmd{int}(M^{m}_{c}-\mathring{M}^{m}_{c})\neq\varnothing$, $\phi_{r}$ has maximal rank $k$ in $\widetilde{U}$ and $k<m$. In analogy to the analysis of maximal rank above, there exists $U\subseteq\widetilde{U}$ such that rank$\phi_{r}\big|_{U}=k$ and $\phi_{r}U=W$ is a $k$-dimensional submanifold. From Lemma \ref{lem7}, $\dim T_{p}W=\dim E_{0}(p)$. There exists a contradiction.
\endproof

\begin{rema}\label{remacon}
From Lemma \ref{lem10} and the constant rank theorem, for any $x\in\mathring{M}^{m}_{c}$, there exists a neighborhood $U(x)$ in $\mathring{M}^{m}_{c}$ such that $W(p)=\phi_{r}U(x)$  is a submanifold of $N$, where $p=\phi_{r}x\in\mathring{M}^{m}_{+}$.
\end{rema}

Due to $H$ has maximal rank $n-m$ in $M^{m}_{+}$, we easily get $M^{m}_{+}$ is open in $M_{+}$. Furthermore, we have
\begin{lemm}\label{lem11}
$\mathring{M}^{m}_{+}$ is dense in $M^{m}_{+}$.
\end{lemm}
\proof Set $U(p)\subseteq M^{m}_{+}=\phi_{r}M^{m}_{c}$ a neighborhood of $\forall p\in M^{m}_{+}$, we only need to prove $U(p)\cap\mathring{M}^{m}_{+}\neq\varnothing$.
$$U(p)\cap\mathring{M}^{m}_{+}=U(p)\cap\phi_{r}\mathring{M}^{m}_{c}\supseteq\phi_{r}(\phi_{r}^{-1}U(p)\cap\mathring{M}^{m}_{c}).$$
From Lemma \ref{lem10}, $\mathring{M}^{m}_{c}$ is dense in $M^{m}_{c}$, $\phi_{r}^{-1}U(p)\cap\mathring{M}^{m}_{c}\neq\varnothing$. Hence, $U(p)\cap\mathring{M}^{m}_{+}\neq\varnothing$.
\endproof

Let $(x^{1}, x^{2}, \dots, x^{n})$ be a local coordinate system at $p$, where $p\in N$. Due to the definition of $H^{*}$ in (\ref{4.4}), we give the expression of $H^{*}$ using local coordinate,
\begin{align}\label{0.01}
H^{*}(\xi)=g^{*ik}(\xi)f_{kj}\xi_{i}dx^{j},
\end{align}
where $\xi\in T_{p}^{*}N$.

\begin{lemm}\label{lem15}
Let $(\xi_{1}, \xi_{2}, \dots, \xi_{n})$ be a local coordinate system of $T_{p}^{*}N$. Then for the linear map ${(dH^{*}_{p})}_{\xi}:T_{\xi}T_{p}^{*}N\rightarrow T_{\xi}T_{p}^{*}N$, we have the following results.\\
(1)For any $\xi\in T^{*}_{p}N$,
\begin{align}\label{4.5}
{(dH^{*}_{p})}_{\xi}=f_{kj}g^{*ik}(p,\xi)\frac{\partial}{\partial \xi_{j}}\otimes d \xi_{i},
\end{align}
(2)For any $p\in\mathring{M}^{m}_{+}$, $\xi\in E^{\bot}_{0}(p)$, $${(dH^{*}_{p})}_{\xi}^{2}=\lambda {(dH^{*}_{p})}_{\xi},\ \ \ \ \ \ {(dH^{*}_{p})}_{\xi}=(n-m)\lambda.$$
\end{lemm}
\proof From (\ref{0.01}), we immediately get the expression of $dH^{*}$ in (1) by directly calculation.

In order to prove (2), for any $p\in\mathring{M}^{m}_{+}$ and $\xi_{0}\in E^{\bot}_{0}(p)=E_{\lambda}^{*}(p)$, let $\xi(t)$ be a curve in $E^{\bot}_{0}(p)$ satisfying $\xi(0)=\xi_{0}$. From Corollary \ref{coro1} and Remark \ref{remacon}, $H^{*}(\xi(t))=\lambda\xi(t)$. Set $X=\xi'(0)\in T_{\xi_{0}}E^{\bot}_{0}(p)$, we have $(dH^{*}_{p})_{\xi_{0}}X=\lambda X$. Define
\begin{align*}
E_{0}^{*}(p,\xi)=\{X^{i}g_{ij}^{*}(p,\xi)\big|X\in E_{0}\}.
\end{align*}
From (\ref{4.5}), $dH^{*}(E_{0}^{*}(p,\xi))=0$, where $\textmd{rank}E_{0}^{*}(\xi)=\textmd{rank}E_{0}=m$.
$$T_{\xi}T^{*}_{p}N=E_{0}^{*}(p,\xi)\oplus T_{\xi}E^{\bot}_{0}(p)$$
and the following algebraic equations hold,
\begin{align}\label{4.6}
{(dH^{*}_{p})}_{\xi}^{2}=\lambda {(dH^{*}_{p})}_{\xi},\ \ \ \ \ \ {(dH^{*}_{p})}_{\xi}=(n-m)\lambda.
\end{align}
This completes the proof.
\endproof

\begin{lemm}\label{lem12}
$M^{m}_{+}$ is open and closed in $M_{+}$, that is, $M^{m}_{+}$ is a union of some connected components of $M_{+}$.
\end{lemm}
\proof
(\ref{4.6}) also holds on the closure of $\mathring{M}^{m}_{+}$ in $M_{+}$ by the continuity of $dH^{*}$. Combine Lemma \ref{lem11}, we have $\overline{\mathring{M}^{m}_{+}}=\overline{M^{m}_{+}}$. Hence, (\ref{4.6}) holds on $\overline{M^{m}_{+}}$ either. Namely, $dH^{*}$ has maximal rank $n-m$ on $\overline{M^{m}_{+}}$. Due to the definition of $M^{m}_{+}$, $M^{m}_{+}$ is closed in $M_{+}$.
\endproof

From Lemma \ref{lem12}, we also have $M^{m}_{c}$ is open and closed in $M_{c}$, that is, $M^{m}_{c}$ is a connected component of $M_{c}$. Replacing $M_{+}$ by $M_{+}\backslash M^{m}_{+}$ and iterating the whole process above by limited steps, we have
\begin{coro}\label{coro2}
(1) $\dim E_{0}(p)$ is constant on each connecting component of $M_{+}$,\\
(2) $\mathring{M}_{+}=\mathop{\cup}\limits_{m\in N^{+}}\mathring{M}^{m}_{+}$ is dense in $M_{+}$.\\
(3) $M_{+}$ is a disjoint union of $M^{m}_{+}$, where $m\in N^{+}$.
\end{coro}

\subsection{Tubes of focal varieties}
\begin{lemm}\label{lem13}
For any $p\in M_{+}^{m}$, we have $\phi_{r}^{-1}(p)=\{\overleftarrow{\exp}_{p}(-r\eta)\big|\eta=\mathcal{L}^{-1}\tau, \tau\in \mathcal{U}E^{\bot}_{0}(p)\}$, which is a $(k-1)$-dimensional closed submanifold in the geodesic sphere $S^{n-1}_{\overleftarrow{F}}(p;r)$, where $k=n-m=\dim E^{\bot}_{0}(p)$.
\end{lemm}
\proof If $p\in\mathring{M}^{m}_{+}$, from Lemma \ref{lem8} and Remark \ref{remacon}, obviously, $\phi_{r}^{-1}(p)=\{\overleftarrow{\exp}_{p}(-r\eta)\big|\eta=\mathcal{L}^{-1}\tau, \tau\in \mathcal{U}E^{\bot}_{0}(p)\}$. Suppose that $p\in M_{+}^{m}\setminus\mathring{M}^{m}_{+}$ and $\tau\in \mathcal{U}E^{\bot}_{0}(p)$. Let $\overleftarrow{\gamma}(s)$ be a geodesic which satisfies $\overleftarrow{\gamma}(0)=p$ and $\overleftarrow{\gamma}'(0)=-\eta=-\mathcal{L}^{-1}(\tau)$. Since $\mathring{M}_{+}^{m}$ is dense in $M_{+}^{m}$ and $E^{\bot}_{0}$ is a smooth subbundle of $T^{*}N$ on $M_{+}^{m}$, there exist $\{p_{i}\}\in \mathring{M}_{+}^{m}$ and $\tau_{i}\in \mathcal{U}E^{\bot}_{0}(p_{i})$ such that $\lim\limits_{i\rightarrow+\infty}p_{i}=p$ and  $\lim\limits_{i\rightarrow+\infty}\tau_{i}=\tau$. $\overleftarrow{\gamma}_{i}(s)$ is the geodesics determined by $(p_{i},-\eta_{i})$ which satisfies $f(\overleftarrow{\gamma_{i}}(r))=c$, where $\eta_{i}=\mathcal{L}^{-1}(\tau_{i})$. Due to continuous dependence of geodesics on initial conditions, $\lim\limits_{i\rightarrow+\infty}\overleftarrow{\gamma}_{i}(s)=\overleftarrow{\gamma}(s)$ for $s\in[0,r]$, $f(\overleftarrow{\gamma}(r))=c$. This implies  $x=\overleftarrow{\gamma}(r)\in M_{c}$. Hence, $\{\overleftarrow{\exp}_{p}(-r\eta)\big|\eta=\mathcal{L}^{-1}\tau, \tau\in \mathcal{U}E^{\bot}_{0}(p)\}\subseteq \phi_{r}^{-1}p$.

Furthermore, there exists a geodesic $\gamma(s)$ determined by $x\in \phi_{r}^{-1}p$ and $\textbf{n}\in\mathcal{N}^{0}_{x}(M_{c})$ which satisfies $\gamma(r)=p$ and $\gamma'(r)=\eta=\mathcal{L}^{-1}(\tau)$. Similar to the proof in (1) of Lemma \ref{lemX}, for any $X\in E_{0}(p)$,
$$0=H^{*}(\tau)(X)=H(\eta, X)=\frac{1}{2}a'(\beta)\tau(X).$$
Namely, we have $\tau\in\mathcal{U}E_{0}^{\bot}(p)$.

Since $\mathcal{L}:SN\rightarrow S^{*}N$ is a diffeomorphism and $\mathcal{U}E_{0}^{\bot}(p)=E_{0}^{\bot}(p)\cap S^{*}_{p}N$, where $E_{0}^{\bot}(p)$ is a $k$-dimensional subspace in $T^{*}_{p}N$, $\{\overleftarrow{\exp}_{p}(-r\eta)\big|\eta=\mathcal{L}^{-1}\tau, \tau\in \mathcal{U}E^{\bot}_{0}(p)\}$ is a closed submanifold in the geodesic sphere $S^{n-1}_{\overleftarrow{F}}(p;r)$.
\endproof

From Lemma \ref{lem13}, the focal map $\phi_{r}:M_{c}\rightarrow M_{+}$ is subjective. Moreover, we have
\begin{lemm}\label{lem0}
For any $x\in M_{c}$, there exist a neighborhood $U(x)$ and a sufficiently small $h\in (0,d(M_{-},M_{+}))$ such that $\gamma_{z}(s)\big|_{[0,r+h]}$ passes $M_{+}$(resp.$M_{-}$) once for any $z\in U(x)$.
\end{lemm}
\proof
For a unique geodesic $\gamma_{x}(s)$ determined by $x=\gamma_{x}(0)$ and $\gamma'_{x}(0)$, set $p=\gamma_{x}(r)\in M_{+}$ and $\eta=\mathcal{L}^{-1}(\tau)=\gamma_{x}'(r)$, where $\tau\in\mathcal{U}E_{0}^{\perp}(p)$. Denote the function $f_{x}(s)=f(\gamma_{x}(s))$, similar to the calculation in Lemma \ref{lem9}, using Taylor expansion for $f_{x}(s)$ at $s=r$,
$$f_{x}(s)=\beta+\frac{1}{2}f_{x}''(\xi)s^2,\ \ \ \ s\in(r,r+d(M_{-},M_{+})),$$

From Lemma \ref{lemX} and Lemma \ref{lem13}, similar to the proof in Corollary \ref{coro1}, $E^{\bot}_{0}(p)=E^{*}_{\lambda}(p)$ for any $p\in M_{+}$. Hence,
$$f''_{x}(r)=H(\eta,\eta)=H^{*}(\tau)(\eta)=\lambda\tau(\eta)=\lambda<0.$$
Then there exists a sufficiently small $h(x)$ such that $f''_{x}(s)<0$ for any $s\in(r,r+h(x))$. Namely, if $\forall s\in(r,r+h(x))$, $f_{x}(s)<\beta$. Since $f_{x}(s)$ is continuous in $x$, for any $x\in M_{c}$, there must exist a neighborhood $U(x)$ of $x$ in $M_{c}$ and a sufficiently small $h\in (0,d(M_{-},M_{+}))$ such that $\gamma_{z}(s)\big|_{[0,r+h]}$ intersect $M_{+}$ once for $z\in U(x)$. This completes the proof.
\endproof
\subsection{Morse Index Theorem}
~~~~Let $(N,F)$ be a Finsler manifold, $P$ a submanifold of $N$ and $\gamma:[0,a]\rightarrow N$ an unit geodesic with $x=\gamma(0)\in P$ and $\textbf{n}(x)=\dot{\gamma}(0)\in \mathcal{N}_{x}^{0}P$. Recall that the Morse index form $I$ along the normal geodesic $\gamma:[0,a]\rightarrow N$ is the symmetric bilinear form
\begin{align}\label{I}
I(X,Y)=\int_{0}^{a}[g_{\dot{\gamma}}(\nabla^{\dot{\gamma}}_{\dot{\gamma}}X,\nabla_{\dot{\gamma}}^{\dot{\gamma}}Y)-g_{\dot{\gamma}}(\textbf{R}_{\dot{\gamma}}(X),Y)]ds,
\end{align}
where $X$, $Y$ are piecewise $C^{\infty}$ vector fields along $\gamma$.

$\mathcal{X}^{P}_{x}$ is the subspace of all piecewise smooth vector fields $X$ along $\gamma$ such that $X(0)\in T_{x}P\setminus\{0\}$, which satisfies $X$ is orthogonal to $\dot{\gamma}$ along $\gamma$ and $X(a)=0$. Denote $I^{P}:\mathcal{X}^{P}_{x}\times\mathcal{X}^{P}_{x}\rightarrow\mathbb{R}$,
\begin{align}
I^{P}(X,Y)=I(X,Y)-g_{\textbf{n}}(A_{\textbf{n}}X(0),Y(0)),
\end{align}
where $A_{\textbf{n}}$ is the shape operator defined in (\ref{A}) with respect to $\textbf{n}$. Define the Morse index $I^{P}$ of $\gamma$ as follows,
$$\textmd{ind}~I^{P}=\max\{\dim(\mathcal{B})~|~\mathcal{B}~\textmd{is}~\textmd{a}~\textmd{subspace}~\textmd{of}~\mathcal{X}^{P}_{x}~\textmd{with}~I^{P}|_{\mathcal{B}}<0\}.$$
A $P$-Jacobi field $J$ is a Jacobi field which satisfies in addition
\begin{align}\label{t}
J(0)\in T_{x}P, \ \ \ \ \nabla_{\textbf{n}}^{\textbf{n}}J(0)+A_{\textbf{n}}J(0)\in \mathcal{N}_{x}^{0}P.
\end{align}
A point $\gamma(s_{0})$, $s_{0}\in[0,a]$ is called a \textit{$P$-focal point} along $\gamma$ if there exists a non-null $P$-Jacobi field $J$ along $\gamma$ with $J(s_{0})=0$. The multiplicity $\mu^{P}(s_{0})$ of a $P$-focal point $\gamma(s_{0})$ is the dimension of the vector space of all $P$-Jacobi fields along $\gamma$ that vanish in $s_{0}$.

Recall the Morse Index Theorem in Finsler geometry,
\begin{lemm}\cite{P}\label{lemXP}
Let $(N,F)$ be a Finsler manifold, $P$ a submanifold of $N$ and $\gamma:[0,a]\rightarrow N$ a geodesic with $\gamma(0)\in P$ and $\dot{\gamma}(0)\in\mathcal{N}^{0}_{x}P$. Then
\begin{align}
\textmd{ind}~I^{P}=\sum\limits_{s_{0}\in(0,a)}\mu^{P}(s_{0})<\infty.
\end{align}
\end{lemm}

\begin{rema}
In \cite{P}, the Morse index form $I$, $P$-Jacobi field $J$ and the the shape operator $A$ are defined in the horizontal subspaces of $T(TN\setminus\{0\})$. By a direct calculation, we know that they can be express on $N$ by (\ref{I}), (\ref{t}) and (\ref{A}).
\end{rema}

Similar to the proof of Lemma 4.6 in $\uppercase\expandafter{\romannumeral2}$ of \cite{Sakai}, we get
\begin{lemm}\label{S}
Let $J(s)$ be a vector field along a normal geodesic $\gamma:[0,a]\rightarrow N$ with $\gamma(0)\in P$ and $\dot{\gamma}(0)\in \mathcal{N}_{x}^{0}P$. Then $J(s)$ is a $P$-Jacobi field if and only if there exists a $C^{\infty}$ variation of $\gamma$ given by
$$\Psi:(-\varepsilon,\varepsilon)\times[0,a]\rightarrow N, \ \ \ \ \Psi(0,s)=\gamma(s),$$
where for each $t\in(-\varepsilon,\varepsilon)$ such that the variation curves $\Psi_{t}(s)=\Psi(t,s)$ are all geodesic perpendicular to $P$ at $s=0$ and $J(s)$ is given by $\frac{\partial\Psi}{\partial t}\big|_{t=0}$.
\end{lemm}

In fact, the $C^{\infty}$ variation in Lemma \ref{S} can be given as $\Psi(t,s)=E(\sigma(t),s\textbf{n}(\sigma(t)))$, where $\sigma(t)$ is a curve in $P$ and $\textbf{n}(\sigma(t))\in\mathcal{N}^{0}_{\sigma(t)}P$. Moreover, from Lemma \ref{lem2} and the proof of Lemma 3.3 in \cite{HCY}, we immediately have
\begin{lemm}\label{J}
For $s_0\in (0,a)$, $p=\gamma(s_{0})$ is a focal point of $P$ if and only if there exists a $P$-Jacobi field along $\gamma(s)$ such that $J(s_{0})=0$.
\end{lemm}
Hence, the focal point defined by the normal exponential map $E$ is equivalent to the $P$-focal point defined by the $P$-Jacobi field.

\section{Proof of Theorem \ref{thm0}}
\begin{lemm}\label{lem21}
The rank of the focal map is constant on each connected component.
\end{lemm}
\proof
Let $\gamma_{x}(s)$ be an $f$-segment through $x\in M_{c}$ which satisfies $\gamma_{x}(0)=x$, $\gamma'(0)=\textbf{n}(x)$, where $\textbf{n}(x)\in\mathcal{N}^{0}_{x}M_{c}$. For any point $z$ which is sufficiently close to $x$, it is obvious that $\textmd{rank}\phi_{r}(z)\geqslant\textmd{rank}\phi_{r}(x)$. In order to prove rank$\phi_{r}$ is constant, we only need to prove rank$\phi_{r}(z)\leq$rank$\phi_{r}(x)$.

From Lemma \ref{lem0}, there exist a neighborhood $U(x)$ of $x$ and a sufficiently small $h\in (0,d(M_{-},M_{+}))$ such that for $z\in U(x)$, the geodesic segment $\Gamma_{z}(s)=\gamma_{z}(s)\big|_{[0,r+h]}$ passes $M_{+}$(resp. $M_{-}$) at $\gamma_{z}(r)$ only. If $\gamma_{x}(r)$ is not a focal point of $M_{c}$ along $\Gamma_{x}(s)$, then rank$\phi_{r}(x)=n-1$ and $n-1=\textmd{rank}\phi_{r}(z)\geq\textmd{rank}\phi_{r}(x)=n-1$. That is, rank$\phi_{r}$ is constant locally. If $\gamma_{x}(r)$ is a focal point of $M_{c}$ along $\Gamma_{x}(s)$, then $\gamma_{x}(r)$ is the unique one. Denote $m(\Gamma_{x}(s))$ as the multiplicity of $\gamma_{x}(r)$, then $m(\Gamma_{x}(s))=n-1-\textmd{rank}_{x}\phi_{r}$. From Lemma \ref{lemXP} and Lemma \ref{J} , $\textmd{ind}I^{M_{c}}_{x}=m(\Gamma_{x}(s))=n-1-\textmd{rank}_{x}\phi_{r}$.

Set $\{e_{i}\}=\{e_{a}, \textbf{n}(x)\}$ an unit orthonormal frame at $x\in M_{c}$ with repect to $g_{\textbf{n}(x)}$. Firstly, translate parallel $\{e_{i}\}$ along $\Gamma_{x}(s)$, the unit orthonormal frame field $\{E_{i}\}=\{E_{a}, T\}$ can be given, where $T=\Gamma'_{x}$. Secondly, translate parallel $\{e_{i}\}$ from $x$ to $z$ along the geodesic on $M_{c}$ with respect to the Riemannian metric $g_{\textbf{n}}$, we can get the unit orthonormal frame $\{\overline{e}_{i}\}=\{\overline{e}_{a}, \textbf{n}(z)\}$ at $z$ with repect to $g_{\textbf{n}(z)}$. Thirdly, translate parallel $\{\overline{e}_{i}\}$ along $\Gamma_{z}(s)$, the unit orthonormal frame field $\{\overline{E}_{i}\}=\{\overline{E}_{a}, \overline{T}\}$ is obtained, where $\overline{T}=\Gamma'_{z}$. Set $X(s)=b^{a}(s)E_{a}(s)$ for any $X\in\mathcal{B}_{x}$, where $b^{a}(r+h)=0$. Set $\overline{X}(s)=b^{a}(s)\overline{E}_{a}(s)\in\mathcal{X}^{M_{c}}_{z}$,
\begin{align}\label{XX}
I^{M{c}}(\overline{X},\overline{X})=\int^{r+h}_{0}\big(\dot{b}_{a}\dot{b}^{a}-b^{a}b^{b}g_{\overline{T}}(R_{\overline{T}}(\overline{E}_{a}),\overline{E}_{b})\big)ds+b^{a}(0)b^{b}(0)g_{\textbf{n}(z)}(\nabla^{\textbf{n}(z)}_{\overline{e}_{a}}\textbf{n}(z),\overline{e}_{b}).
\end{align}
Since $\gamma_{z}(s)$, $\overline{E}_{a}$ and $\overline{T}$ are determined by $z\in M_{c}$, from (\ref{XX}), $I(\overline{X},\overline{X})$ only depends on $z$ continuously.
Due to $I(X,X)<0$, there exists $\overline{U}(x)\subseteq U(x)$ such that $I(\overline{X},\overline{X})<0$ for all $z\in\overline{U}(x)$. Namely, $\overline{X}\in\mathcal{B}_{z}$ if $z\in\overline{U}(x)$. Hence, $\textmd{dim}\mathcal{B}_{z}\geq\textmd{dim}\mathcal{B}_{x}$. Therefore, $$m(\Gamma_{z})=\textmd{ind}I^{M_{c}}_{z}=\max\textmd{dim}\mathcal{B}_{z}\geq\max\textmd{dim}\mathcal{B}_{x}=\textmd{ind}I^{M_{c}}_{x}=m(\Gamma_{x}).$$
Due to the expression of $m(\Gamma_{x})$, rank$\phi_{r}(z)\leq$rank$\phi_{r}(x)$. Hence, The rank of the focal map is locally constant.

Furthermore, if rank$\phi_{r}\big|_{U(x)}=m$, then we have $\phi_{r}\big|_{\overline{U(x)}}=m$. Hence, the set of $\phi_{r}$ with constant rank is open and closed. Namely, it is a connected component of $M_{c}$. This completes the proof.
\endproof

Combine Lemma \ref{lem9} and Lemma \ref{lem21}, $M^{m}_{+}$ is a submanifold of $(N,F)$. According to the analysis of Lemma \ref{lem8} about tubes, $\phi_{r}^{-1}M^{m}_{+}$ is a tube of radius $r$ with respect to $\overleftarrow{F}$. We completes the proof of Theorem \ref{thm0} while the focal variety is $M_{+}$. When the focal variety is $M_{-}$, we just need to change $F$ to $\overleftarrow{F}$ in all proofs above.

From Theorem \ref{thm0}, similar to the proof in \cite{GT1}, we have
\begin{coro}
Each connected component of $M_{\pm}$ has codimension not less than 2 if and only if $M_{\pm}=\textmd{Foc}(f)$.
\end{coro}

\section{Example}
~~~~Recall $(\alpha,\beta)$-metric in Finsler geometry which can be expressed by $F=\alpha\phi(s)$, where $\alpha=\sqrt{a_{ij}(x)y^{i}y^{j}}$ is a Riemannian metric, $\beta=b_{i}(x)y^{i}$ is a 1-form satisfying $b:=\|\beta\|_{\alpha}<b_{0}$ $s=\frac{\beta}{\alpha}$ for any $x\in N$ and $\phi:(-b_{0},b_{0})\rightarrow(0,+\infty)$ is a smooth positive function, which satisfies
\begin{align}\label{a}
\phi(s)>0,\ \ \ \ (\phi(s)-s\phi'(s))+(b^2-s^2)\phi''(s)>0,
\end{align}
where $s$ and $b$ are arbitrary numbers with $|s|\leq b<b_{0}$.

Let $\alpha$ be a standard Euclidean sphere metric, which its dual metric denotes by $\alpha^{*}$. Set $\varphi(s):[-1,1]\rightarrow(0,+\infty)$ a positive function which is expressed by $\varphi(s)=1+\varepsilon s^{k}$, where $0<\varepsilon<\frac{1}{2k-1}$ and $k\in N_{+}$. By a direct computation, we can verify that (\ref{a}) holds. Denote $x=(\overline{x}_{1}, \overline{x}_{2})\in \mathbb{R}^{n+1}$, where $\overline{x}_{1}=(x_{1}, x_{2},\dots, x_{p+1})$, $\overline{x}_{2}=(x_{p+2}, x_{p+3}\dots, x_{p+q+2})$, $p,q\in N^{+}$ and $p+q=n-1$. For $x\in\mathbb{S}^{n}$, define $V(x)$ as the projection of $\widetilde{x}=(\overline{x}_{1}, -\overline{x}_{2})$ onto $\mathbb{S}^{n}$, then $V(x)=\tilde{x}-\langle\tilde{x},x\rangle x$. Set $\beta^{*}(\xi)=\xi(V)$. Then we equip an $(\alpha,\beta)$-metric $F$ on $\mathbb{S}^{n}$, which its dual metric is expressed by $F^{*}=\alpha^{*}\varphi(s)$, where $s=\frac{\beta^{*}}{\alpha^{*}}$.

Set $\Phi(x)=|\overline{x}_{1}|^2-|\overline{x}_{2}|^{2}$ and $f=\Phi\big|_{\mathbb{S}^{n}}$, we have
$$\big(\alpha^{*}(df)\big)^2=|\nabla^{\alpha}f|^{2}=\langle \nabla^{E}\Phi-2\Phi x, \nabla^{E}\Phi-2\Phi x\rangle=4(1-f^{2}),$$
$$\beta^{*}(df)=\langle \nabla^{\alpha}f,V\rangle_{\alpha}=\langle \nabla^{E}\Phi-2\Phi x,\tilde{x}-\langle\tilde{x},x\rangle x\rangle=2(1-f^2).$$
Hence,
$$F(\nabla f)=F^{*}(df)=\alpha^{*}(df)\varphi(\frac{\beta^{*}(df)}{\alpha^{*}(df)})=2\sqrt{1-f^{2}}(1+\varepsilon(1-f^2)^{\frac{k}{2}}).$$
Namely, $f$ is a $C^{\infty}$ transnormal function on $(\mathbb{S}^{n},F)$. In this case, $J=[-1,1]$.
\begin{exam}\label{exam1}
If $k=1$,
$$a(t)=4(1-f^2)+8\varepsilon(1-f^{2})^{\frac{3}{2}}+4\varepsilon^{2}(1-f^{2})^{2}\in C^{1}(J).$$
Hence, $a(t)\neq0$ for $t\in(-1,1)$ and $a'(\pm1)=\mp8\neq0$.
In this case, $f$ and $a(t)$ satisfy the conditions in Theorem \ref{thm0}.
\end{exam}
\begin{exam}\label{exam2}
If $k=3$,
$$a(t)=4(1-f^2)+8\varepsilon(1-f^{2})^{\frac{5}{2}}+4\varepsilon^{2}(1-f^{2})^{4}\in C^{2}(J),$$
If $k=2l$, $l\in N_{+}$,
$$a(t)=4(1-f^2)+8\varepsilon(1-f^{2})^{l+1}+4\varepsilon^{2}(1-f^{2})^{2l+1}\in C^{\infty}(J).$$
Hence, there exists a smooth function $f$ on $N=\mathbb{S}^{n}$ which satisfies $a\in C^{2}(J)$.
\end{exam}

Yali Chen \\
School of Mathematical Sciences, Tongji University, Shanghai, 200092, China\\
E-mail: chenylwuhu@qq.com\\

Qun He \\
School of Mathematical Sciences, Tongji University, Shanghai, 200092, China\\
E-mail: hequn@tongji.edu.cn\\

\end{document}